\documentclass[11pt,a4paper]{article}
\usepackage[T1]{fontenc}
\usepackage[margin=27mm]{geometry}
\usepackage{amsmath,amssymb,amsthm,mathtools}
\usepackage{enumitem}
\usepackage{tikz}
\usepackage[hidelinks]{hyperref}
\usepackage{fancyhdr}
\hypersetup{pdftitle={Reflection positivity of discrete and Brownian loop-soup occupation fields}}
\setlist[enumerate]{itemsep=2pt,topsep=4pt}
\numberwithin{equation}{section}
\newtheorem{theorem}{Theorem}[section]
\newtheorem{lemma}[theorem]{Lemma}
\newtheorem{proposition}[theorem]{Proposition}

\theoremstyle{remark}
\newtheorem{remark}[theorem]{Remark}
\newcommand{\R}{\mathbb R}
\newcommand{\E}{\mathbb E}
\newcommand{\Sym}{\operatorname{Sym}}
\newcommand{\SPD}{\mathcal P}
\newcommand{\diag}{\operatorname{diag}}
\newcommand{\tr}{\operatorname{tr}}
\newcommand{\rank}{\operatorname{rank}}
\newcommand{\T}{\mathsf T}
\newcommand{\W}{\mathcal W}
\newcommand{\dd}{\mathrm d}
\title{Reflection positivity of discrete and Brownian loop-soup occupation fields}
\author{Yves Le Jan\\ Université Paris-Saclay}
\date{24 November 2026}

\begin{document}
\maketitle
\thispagestyle{plain}

\begin{abstract}
The intensities guaranteeing reflection positivity of loop-soup
vertex occupations on every reflection-symmetric finite graph with
positive semidefinite cross-conductance matrix of rank $r$ form the
real Wallach set $\W_r$.
For the renormalized occupation field on the massless planar Dirichlet
strip, the exact set is
$\bigcap_{r\geq1}\W_r=\{0,\tfrac12,1,\tfrac32,\ldots\}$.
\end{abstract}

\section{Introduction}

Reflection positivity is a necessary condition for a Euclidean
statistical field to admit a unitary quantum interpretation by Wick
rotation. In the reconstruction of Osterwalder and Schrader, the
expected product of an observable supported at positive times with its
reflected complex conjugate gives the squared norm of the corresponding
quantum state. Under the other reconstruction hypotheses, Euclidean time
translations act on the resulting Hilbert space by $e^{-tH}$, where
$H$ is a nonnegative self-adjoint Hamiltonian. Wick rotation gives the
unitary evolution $e^{-isH}$ in real time $s$
\cite{OsterwalderSchraderI,OsterwalderSchraderII}. A negative value of
the reflection form therefore rules out this interpretation for the
given field and choice of time reflection.

The author's 2006 isomorphism identifies the vertex occupation field $L_1$
of a Markov loop soup of intensity one on a finite graph with
$\tfrac12|\Phi|^2$ in law, where $\Phi=\phi_1+i\phi_2$ for two
independent real Gaussian free fields \cite[Proposition~2]{LeJan2006}. In its real form, the
occupation field $L_{1/2}$ at intensity $1/2$ has the law of
$\phi^2/2$ for a real Gaussian free field $\phi$
(see \cite[Theorem~2(a)]{LeJan} or \cite[Section~6.2]{LeJanb}).
By superposition, $L_{N/2}$ has the law of $\tfrac12\sum_{j=1}^N\phi_j^2$
for independent copies $\phi_j$ of $\phi$. Suppose a reflection
exchanges the two halves of the graph, preserves conductances and
killing, and has a positive semidefinite cross-conductance matrix.
Standard Gaussian reflection positivity, together with the preceding
isomorphism, gives reflection positivity at every intensity
$N/2$, $N=1,2,\ldots$ \cite[Remark~25(c)]{LeJan}. We shall use this
standard sufficiency statement without further discussion.

On a subdivided cube the author found a function of the local times of a
single loop whose reflection form under the loop measure is strictly
negative. On this graph the occupation field fails reflection positivity
at sufficiently small positive intensities
(see \cite[Section~9.2]{LeJan} or \cite[Section~12.2]{LeJanb}). Which positive intensities guarantee
reflection positivity on every graph satisfying the conductance
condition? Theorem~\ref{thm:main} shows that half integers are the only such intensities.
At fixed cross-conductance rank $r\geq1$, every intensity above
$(r-1)/2$ also gives reflection positivity. Together with the
nonnegative multiples of $1/2$, these intensities form the real Wallach
set $\W_r$. For each $r$ we construct a graph whose admissible set is
exactly $\W_r$.

The essential new ingredient is the use of the classical matrix gamma
formula \eqref{eq:gamma} for $\det\!\left(X+Y\right)^{-\alpha}$.
Testing the occupations $L_+$ and $L_-$ on the two halves by nonnegative
potentials $s$ and $t$ gives a joint Laplace transform of the form
\[
 \E_\alpha e^{-\langle s,L_+\rangle-\langle t,L_-\rangle}
       =\frac{b_\alpha(s)b_\alpha(t)}{\det(X_s+X_t)^\alpha},
\]
where $b_\alpha$ is positive and $X_s,X_t$ are positive definite
symmetric matrices of order $r$. For $\alpha>(r-1)/2$, the gamma
integral gives reflection positivity by expressing the determinant
factor as the $L^2$ inner product of $T\mapsto e^{-\tr(X_sT)}$ and
$T\mapsto e^{-\tr(X_tT)}$.
On the graphs we construct, varying the vertex potentials makes $X_s$
range over an open set. After a change of matrix coordinates and
restriction to an $m\times m$ block, reflection positivity would
therefore imply that the kernel $\det(X+Y)^{-\alpha}$ is positive
definite near $X=Y=I_m/2$. To test this, we apply the same differential
operator to each argument. In the $X$ variable this operator is
$D_X=\det(\partial_X^{\mathrm{sym}})$, the determinant of the symmetric
matrix of partial derivatives, with the off-diagonal derivatives
multiplied by $1/2$. Positive definiteness requires
\[
 \left.D_XD_Y\,\det(X+Y)^{-\alpha}\right|_{X=Y=I_m/2}\geq 0,
\]
as one sees by approximating the operator on each side by the same
finite differences. Lemma~\ref{lem:gamma-differential} evaluates the
left-hand side using the matrix gamma integral and a polynomial
identity. For $\alpha\notin\W_r$ its value is negative when
$m=\lfloor 2\alpha\rfloor+2\leq r$.

On the massless planar strip with Dirichlet boundary conditions, let
$X_\alpha$ denote the centered occupation field of the Brownian loop
soup of intensity $\alpha$. (Recall that to define $X_\alpha(f)$ for a smooth,
compactly supported test function $f$, one sums the occupation integrals
of $f$ over loops of duration greater than $\delta$, subtracts the mean
of this sum, and takes the $L^2$ limit as $\delta\downarrow0$.)
For each finite $m$, we recover the mixed determinant derivative above
as a limit involving finite linear combinations of
$e^{-X_\alpha(f)}$, with $f$ nonnegative and supported in the positive
half-strip. Reflection positivity would therefore require this derivative
to be nonnegative for every $m$, which rules out every intensity for
which $2\alpha$ is not an integer.

\subsection*{Finite graphs and the main theorem}

Let $V$ be a finite vertex set with symmetric conductances
$C_{xy}=C_{yx}\geq0$ for $x\ne y$, $C_{xx}=0$, and killing rates
$\kappa_x\geq0$. Write
\[
 Q_{xx}=\kappa_x+\sum_{y\ne x}C_{xy},
 \qquad Q_{xy}=-C_{xy}\quad(x\ne y),
\]
and assume $Q\succ0$. Consider the continuous-time Markov chain with
generator $-Q$, including killing. The distribution of the total vertex local times
$L_\alpha=(L_\alpha^x)_{x\in V}$ of its loop soup of intensity $\alpha>0$
is given by the Laplace identity
\begin{equation}\label{eq:laplace}
 \E_\alpha e^{-\langle s,L_\alpha\rangle}
 =\left(\frac{\det Q}{\det(Q+\diag s)}\right)^\alpha,
 \qquad s\in[0,\infty)^V.
\end{equation}
We set $L_0=0$. Identity \eqref{eq:laplace} is standard for Markovian
loop soups; see \cite[Chapter~4]{LeJan}. 
%A reflection-invariant positive rescaling of the local times does not change reflection positivity.

Suppose that $V=V_+\sqcup V_-$ and that an involution $\rho$ exchanges
$V_+$ and $V_-$, preserves the conductances and killing, and satisfies
\begin{equation}\label{eq:cross}
 B=(C_{x,\rho(y)})_{x,y\in V_+}\succeq0.
\end{equation}
After identifying $V_-$ with $V_+$ through $\rho$, we have
\begin{equation}\label{eq:block}
 Q=\begin{pmatrix}A&-B\\-B&A\end{pmatrix}.
\end{equation}
We have $A-B\succ0$. Indeed, let $n=|V_+|$ and take any
$v\in\R^n\setminus\{0\}$. Since $Q\succ0$,
\[
 0<
 \begin{pmatrix}v\\v\end{pmatrix}^{\T}
 Q
 \begin{pmatrix}v\\v\end{pmatrix}
 =2v^{\T}(A-B)v.
\]
Let $L_+=(L_\alpha^x)_{x\in V_+}$ and
$L_-=(L_\alpha^{\rho(x)})_{x\in V_+}$. The occupation field is
\emph{reflection positive} if
\begin{equation}\label{eq:rp}
 \E_\alpha\big[\overline{F(L_+)}F(L_-)\big]\geq0
 \qquad\text{for every }F\in L^2(\operatorname{law}(L_+)).
\end{equation}
The reflection symmetry gives the same marginal law for $L_+$ and $L_-$.
The finite-graph involutions in Theorem~\ref{thm:main} have no fixed vertices.

For $r\geq1$, define the real Wallach set \cite{Gindikin,WallachII}
\begin{equation}\label{eq:wallach}
 \W_r=
 \left\{\frac j2: j=0,1,\ldots,r-1\right\}
 \cup\left(\frac{r-1}{2},\infty\right),
 \qquad \W_0=[0,\infty).
\end{equation}

\begin{theorem}\label{thm:main}
The following statements hold.
\begin{enumerate}[label=\textup{(\roman*)}]
\item For every graph satisfying \eqref{eq:cross}, the vertex occupation
field is reflection positive whenever $\alpha\in\W_r$, where
$r=\rank B$.
\item For every integer $r\geq1$ and every $\kappa>0$, there is a connected
finite graph $G_r$ with unit edge conductances, killing $\kappa$ at every
vertex, and $\rank B=r$, for which reflection positivity holds if and only
if $\alpha\in\W_r$.
\item A positive intensity gives reflection positivity on every finite
graph satisfying \eqref{eq:cross} if and only if it belongs to
$\{\tfrac12,1,\tfrac32,\ldots\}$.
\end{enumerate}
\end{theorem}

\section{Laplace test functions and the matrix gamma integral}
\label{sec:factorization}

\subsection{The determinant factorization}

For $s\in[0,\infty)^{V_+}$ put
$f_s(\ell)=e^{-\langle s,\ell\rangle}$, and define
\begin{equation}\label{eq:Kdef}
 K_\alpha(s,t)=\E_\alpha[f_s(L_+)f_t(L_-)].
\end{equation}
A kernel is called positive definite when all its finite Gram matrices
are positive semidefinite.

\begin{lemma}\label{lem:laplace}
The occupation field is reflection positive if and only if
$K_\alpha$ is a positive definite kernel on $[0,\infty)^{V_+}$.
\end{lemma}

\begin{proof}
For $F=\sum_i c_i f_{s_i}$ the reflection form is
$\sum_{i,j}\overline{c_i}c_jK_\alpha(s_i,s_j)$. The span of the
$f_s$ is dense in $L^2(\operatorname{law}(L_+))$. Indeed, under
$\ell\mapsto(e^{-\ell_x})_{x\in V_+}$, the functions with
$s\in\{0,1,\ldots\}^{V_+}$ become the coordinate monomials.
Their span consists of the polynomials, which are dense in $L^2$
of the image measure on $[0,1]^{V_+}$. Finally, equality of the two marginals and
Cauchy--Schwarz give
\[
 \left|\E_\alpha[\overline{F(L_+)}G(L_-)]\right|
 \leq\|F\|_2\|G\|_2,
\]
so positivity on this span extends to the full space.
\end{proof}

Let $n=|V_+|$ and $r=\rank B>0$. The spectral theorem gives
\[
 B=\sum_{j=1}^r\lambda_j u_ju_j^{\T}=RR^{\T},\qquad
 R=(\sqrt{\lambda_1}u_1\ \cdots\ \sqrt{\lambda_r}u_r)
       \in\R^{n\times r},
\]
where $\lambda_j>0$ and the $u_j$ are orthonormal. In particular,
the columns of $R$ are linearly independent. Put
\begin{equation}\label{eq:PH}
 P_s=A-B+\diag s,\qquad
 H_s=R^{\T}P_s^{-1}R,\qquad X_s=\tfrac12 I_r+H_s.
\end{equation}
All these matrices are positive definite.

By \eqref{eq:laplace} the joint Laplace transform is
$K_\alpha(s,t)=(\det Q/\det Q_{s,t})^\alpha$, where
\[
 Q_{s,t}=\begin{pmatrix}A+\diag s&-B\\-B&A+\diag t\end{pmatrix}
        =D_{s,t}+UU^{\T},\qquad
 D_{s,t}=\diag(P_s,P_t),\qquad
 U=\begin{pmatrix}R\\-R\end{pmatrix}.
\]
Factoring out $D_{s,t}$ and applying Sylvester's determinant identity
$\det(I+CD)=\det(I+DC)$ gives
\[
\begin{aligned}
 \det Q_{s,t}
 &=\det D_{s,t}\,
       \det(I_{2n}+D_{s,t}^{-1}UU^{\T})\\
 &=\det P_s\det P_t\,
       \det(I_r+U^{\T}D_{s,t}^{-1}U)\\
 &=\det P_s\det P_t\det(I_r+H_s+H_t)\\
 &=\det P_s\det P_t\det(X_s+X_t).
\end{aligned}
\]
Therefore
\begin{equation}\label{eq:factor}
 K_\alpha(s,t)
 =\det(Q)^\alpha\det(P_s)^{-\alpha}\det(P_t)^{-\alpha}
                      \det(X_s+X_t)^{-\alpha}.
\end{equation}
Multiplication by positive factors depending separately on $s$ and $t$
preserves positive definiteness.

\subsection{The matrix gamma integral and sufficiency}

Let $\SPD_r$ be the cone of real symmetric positive definite
$r\times r$ matrices. For $a>(r-1)/2$, define
\begin{equation}\label{eq:matrix-gamma}
 \Gamma_r(a)=\int_{\SPD_r}e^{-\tr U}
          (\det U)^{a-(r+1)/2}\,\dd U
 =\pi^{r(r-1)/4}\prod_{j=0}^{r-1}\Gamma(a-j/2),
\end{equation}
where $\dd U=\prod_{i\leq j}\dd U_{ij}$. The product follows by
the Cholesky substitution $U=LL^{\T}$ with $L$ lower triangular
and $L_{ii}>0$. Its Jacobian is $2^r\prod_iL_{ii}^{r+1-i}$.
The off-diagonal Gaussian integrals give $\pi^{r(r-1)/4}$,
and the integral in $L_{ii}$ contributes
$\tfrac12\Gamma(a-(i-1)/2)$, cancelling one factor of $2$.

For $S=X+Y\succ0$, substitute $U=S^{1/2}TS^{1/2}$. Then
\[
 \tr U=\tr(ST),\qquad \det U=\det S\det T,\qquad
 \dd U=(\det S)^{(r+1)/2}\,\dd T.
\]
For the last identity, diagonalize $S^{1/2}$, with eigenvalues
$a_1,\ldots,a_r$. The independent symmetric coordinates are
multiplied by $a_i^2$ on the diagonal and $a_ia_j$ off the
diagonal. The product of these factors is $\prod_i a_i^{r+1}$.
Substitution into the defining integral gives
\[
 \Gamma_r(a)=(\det S)^a\int_{\SPD_r}e^{-\tr(ST)}
                    (\det T)^{a-(r+1)/2}\,\dd T.
\]
Thus
\begin{equation}\label{eq:gamma}
 \det(X+Y)^{-a}
 =\frac1{\Gamma_r(a)}\int_{\SPD_r}
       e^{-\tr(XT)}e^{-\tr(YT)}
       (\det T)^{a-(r+1)/2}\,\dd T.
\end{equation}

\begin{proof}[Proof of Theorem~\ref{thm:main}\textup{(i)}]
The case $\alpha=0$ follows from $L_0=0$. If $r=0$, then $B=0$
and the joint Laplace transform factors. The two halves are thus
independent and identically distributed, so the reflection form is
$|\E_\alpha F(L_+)|^2$. For $r>0$, the positive half-integer
cases are proved in \cite[Remark~25(c)]{LeJan} (see also
\cite[Section~12.1]{LeJanb}). Every other
$\alpha\in\W_r$ satisfies $\alpha>(r-1)/2$, so
\eqref{eq:gamma} and \eqref{eq:factor} make $K_\alpha$ positive
definite. Lemma~\ref{lem:laplace} concludes the proof.
\end{proof}

\section{The determinant-kernel criterion}
Throughout this section, $\Sym_r$ denotes the vector space of real
symmetric $r\times r$ matrices, and $\SPD_r\subset\Sym_r$ its cone
of positive-definite matrices. We use the independent coordinates
$(x_{ij})_{1\le i\le j\le r}$, with $x_{ji}=x_{ij}$.
%The parameter set occurs in the work of Gindikin and Wallach
%\cite{Gindikin,WallachII}. Sra states the whole-cone matrix-kernel
%criterion \cite[Theorem~3.10]{Sra}. The local criterion also
%follows by restricting to an open convex neighbourhood and applying
%the Laplace representation in Scott and Sokal
%\cite[Theorem~6.4]{ScottSokal}. Uniqueness of the Laplace transform
%identifies the representing measure with the Riesz distribution for
%the inverse determinant power, whose positivity is characterized in
%\cite[Theorem~4.8]{ScottSokal}. We give a direct proof using the
%matrix gamma integral.
%
For every real $a$, define
\begin{equation}\label{eq:gamma-polynomial}
 p_r(a)=\prod_{j=0}^{r-1}(a-j/2).
\end{equation}
The product formula for $\Gamma_r$ gives, for $a>(r-1)/2$,
\begin{equation}\label{eq:gamma-ratios}
 \frac{\Gamma_r(a+1)}{\Gamma_r(a)}=p_r(a).
\end{equation}

For symmetric matrix coordinates, put
\[
 (\partial_X^{\mathrm{sym}})_{ij}
 =
 \begin{cases}
  \partial/\partial x_{ii},&i=j,\\
  \tfrac12\partial/\partial x_{ij},&i<j,\\
  \tfrac12\partial/\partial x_{ji},&i>j,
 \end{cases}
 \qquad D_X=\det(\partial_X^{\mathrm{sym}}).
\]
Thus $D_X$ is a real constant-coefficient differential operator.
For instance, in dimension two,
\[
 D_X=\partial_{11}\partial_{22}
                  -\tfrac14\partial_{12}^2.
\]
For symmetric matrices $X$ and $T$,
\[
 \tr(XT)=\sum_i x_{ii}t_{ii}
              +2\sum_{i<j}x_{ij}t_{ij}.
\]
Thus differentiation with respect to an off-diagonal coordinate
$x_{ij}$ produces the factor $-2t_{ij}$. The normalization by
$1/2$ ensures that, entrywise,
\[
 \partial_X^{\mathrm{sym}}e^{-\tr(XT)}
       =-T e^{-\tr(XT)}.
\]
Since the coordinate derivatives commute, taking the determinant
of the differential-operator matrix gives
\[
 D_Xe^{-\tr(XT)}
       =(-1)^r\det(T)e^{-\tr(XT)}.
\]

\begin{lemma}\label{lem:gamma-differential}
For every real $a$,
\begin{equation}\label{eq:gamma-differential}
 \left.D_XD_Y\det(X+Y)^{-a}\right|_{X=Y=I_r/2}
       =p_r(a)p_r(a+1).
\end{equation}
\end{lemma}

\begin{proof}
For $a>(r-1)/2$ we may differentiate \eqref{eq:gamma} under
the integral. For $X+Y\succeq I_r/2$, derivatives through order
$2r$ introduce polynomial factors bounded by a constant times
$(1+\tr T)^{2r}$. The remaining exponential $e^{-\tr(T)/2}$
absorbs this polynomial, and the resulting bound is integrable
by \eqref{eq:matrix-gamma} after rescaling. Hence
\[
 \begin{aligned}
 &\left.D_XD_Y\det(X+Y)^{-a}\right|_{X=Y=I_r/2}\\
 &\qquad=
 \frac1{\Gamma_r(a)}
 \int_{\SPD_r}e^{-\tr T}
       (\det T)^{a+2-(r+1)/2}\,\dd T\\
 &\qquad=\frac{\Gamma_r(a+2)}{\Gamma_r(a)}
        =p_r(a)p_r(a+1).
 \end{aligned}
\]
The left-hand side is a polynomial in $a$, as follows by
differentiating $\exp(-a\log\det(X+Y))$ and using
$\det(X+Y)=1$ at the evaluation point.
Equality of polynomials therefore proves
\eqref{eq:gamma-differential} for every real $a$.
\end{proof}

Sra gives the whole-cone criterion in matrix form
\cite[Theorem~3.10]{Sra}. On an open convex neighbourhood, uniqueness
of the Laplace transform identifies the measure in Scott and Sokal's
representation theorem with the Riesz distribution. Its positivity
gives the same parameter set \cite[Theorems~6.4 and~4.8]{ScottSokal}.
We give a direct proof using Lemma~\ref{lem:gamma-differential}.

\begin{proposition}\label{prop:local}
For $\alpha\geq0$, the kernel
\[
 J_\alpha(X,Y)=\det(X+Y)^{-\alpha}
\]
is positive definite on a nonempty open subset of $\SPD_r$
if and only if $\alpha\in\W_r$. In that case it is positive
definite on the whole cone.
\end{proposition}

\begin{proof}
Sufficiency follows from the Gaussian integral for nonnegative
half-integers and from \eqref{eq:gamma} for
$\alpha>(r-1)/2$.

We first record the consequence of positive definiteness that
will be used for necessity. Let $k$ be a smooth positive
definite kernel near $(x_0,x_0)$, and let $D$ be any real
constant-coefficient differential operator. Replacing its
coordinate derivatives by finite differences gives an approximation of $D$ of the form
\[
 D_hf(x_0)=\sum_\ell b_\ell(h)f(x_\ell(h)),
 \qquad b_\ell(h)\in\R,\qquad x_\ell(h)\longrightarrow x_0.
\]
Positive definiteness implies
\[
 (D_{h,x}D_{h,y}k)(x_0,x_0)
 =\sum_{\ell,j}b_\ell(h)b_j(h)
                  k(x_\ell(h),x_j(h))\geq0.
\]
Taking $h\downarrow0$ yields
\begin{equation}\label{eq:differential-positivity}
 (D_xD_yk)(x_0,x_0)\geq0.
\end{equation}
This applies in particular to the determinant differential
operator defined above.

%It is useful to see this implication in the simplest case. For a
%coordinate direction $v$, take
%\[
% D_hf(x_0)=\frac{f(x_0+hv)-f(x_0)}{h}.
%\]
%Then positive definiteness gives
%\[
%\begin{aligned}
% 0\leq{}&\frac1{h^2}\bigl[
% k(x_0+hv,x_0+hv)-k(x_0+hv,x_0)\\
% &\hspace{34mm}-k(x_0,x_0+hv)+k(x_0,x_0)
% \bigr],
%\end{aligned}
%\]
%and the limit is $(\partial_v^x\partial_v^yk)(x_0,x_0)$. A higher
%order operator is handled in exactly the same way: write a finite
%difference approximation as a signed linear combination of point
%evaluations, and use the same coefficients in the two kernel variables.
%The signs inside the operator cause no difficulty because the resulting
%expression is still a quadratic form of a positive semidefinite Gram
%matrix. For example, when $r=2$ this applies to
%\[
% D_X=\partial_{11}\partial_{22}-\tfrac14\partial_{12}^2.
%\]
%
Now suppose $J_\alpha$ is positive definite on an open set
containing $X_0\succ0$. With $C=(2X_0)^{1/2}$, the substitution
$X=C\widetilde X C$, $Y=C\widetilde Y C$ moves the center to
$I_r/2$ and multiplies the kernel by a positive constant.

If $\alpha\notin\W_r$, put
\[
 m=\lfloor2\alpha\rfloor+2.
\]
Then
\[
 2\leq m\leq r,\qquad
 \frac{m-2}{2}<\alpha<\frac{m-1}{2}.
\]
Restrict the kernel to matrices of the form
\[
 \widetilde X=\diag(X',I_{r-m}/2),\qquad
 \widetilde Y=\diag(Y',I_{r-m}/2),
\]
where $X',Y'\in\SPD_m$ vary in a sufficiently small neighbourhood
of $I_m/2$. 
%The second block is absent when $m=r$.
This restriction gives the positive
definite kernel $\det(X'+Y')^{-\alpha}$ near $I_m/2$.

Lemma~\ref{lem:gamma-differential}, in dimension $m$, gives
\[
 \left.D_{X'}D_{Y'}\det(X'+Y')^{-\alpha}
          \right|_{X'=Y'=I_m/2}
 =p_m(\alpha)p_m(\alpha+1)<0.
\]
Indeed, exactly the last factor of $p_m(\alpha)$ is negative,
whereas every factor of $p_m(\alpha+1)$ is positive.
This contradicts \eqref{eq:differential-positivity}.

For the first two forbidden intervals, when the ambient rank contains
the indicated block, the obstruction is especially transparent. If
$r\geq2$ and $0<\alpha<1/2$, take $m=2$ and obtain
\[
 p_2(\alpha)p_2(\alpha+1)
 =\alpha\left(\alpha-\tfrac12\right)
  (\alpha+1)\left(\alpha+\tfrac12\right)<0.
\]
If $r\geq3$ and $1/2<\alpha<1$, take $m=3$; then
\[
 p_3(\alpha)p_3(\alpha+1)
 =\alpha^2\left(\alpha-\tfrac12\right)(\alpha-1)
  \left(\alpha+\tfrac12\right)(\alpha+1)<0.
\]
The choice $m=\lfloor2\alpha\rfloor+2$ is the same sign argument on
the unique half-integer interval containing a general
$\alpha\notin\W_r$.
\end{proof}

\section{Graphs realizing the obstruction}
\subsection{The graphs \texorpdfstring{$G_r$}{G\_r}}
Fix $r\geq1$ and $\kappa>0$. Construct the positive half of $G_r$ from
vertices
\[
 v_1,\ldots,v_r,\qquad
 u_{ij}\quad(1\leq i<j\leq r).
\]
Join $u_{ij}$ to $v_i$ and $v_j$. Take a reflected copy, and add the
crossing edges $v_i\leftrightarrow\rho(v_i)$ for $1\leq i\leq r$.
There are no other edges. The vertices $v_1,\ldots,v_r$
form the interface in $V_+$, since they are precisely the vertices
incident to crossing edges.

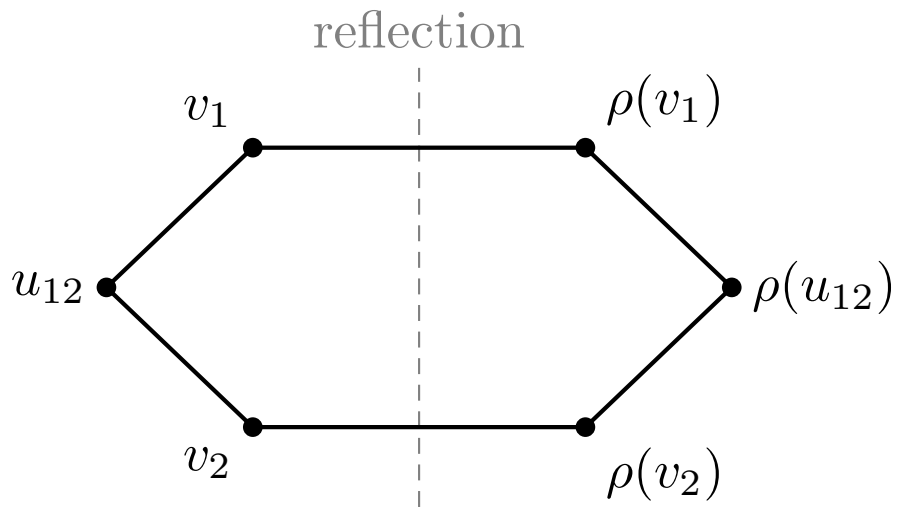
\begin{figure}[htbp]
\centering
\begin{tikzpicture}[scale=0.88,every node/.style={font=\small}]
 \draw[densely dashed,gray] (0,-1.65)--(0,1.65);
 \node[gray,above] at (0,1.65) {reflection};
 \coordinate (a) at (-1.25,1.05);
 \coordinate (b) at (-2.35,0);
 \coordinate (c) at (-1.25,-1.05);
 \coordinate (d) at (1.25,1.05);
 \coordinate (e) at (2.35,0);
 \coordinate (f) at (1.25,-1.05);
 \draw[thick] (a)--(b)--(c)--(f)--(e)--(d)--(a);
 \foreach \p in {a,b,c,d,e,f} \fill (\p) circle (2.1pt);
 \node[above left] at (a) {$v_1$};
 \node[left] at (b) {$u_{12}$};
 \node[below left] at (c) {$v_2$};
 \node[above right] at (d) {$\rho(v_1)$};
 \node[right] at (e) {$\rho(u_{12})$};
 \node[below right] at (f) {$\rho(v_2)$};
\end{tikzpicture}
\caption{The rank-two example $G_2$. All conductances are one and
all killing rates are $\kappa>0$.}
\label{fig:sixcycle}
\end{figure}

Give every edge conductance one and every
vertex killing $\kappa$. This graph is connected, has $r(r+1)$
vertices, and its reflection exchanges the two copies. With the
interface vertices listed first,
\begin{equation}\label{eq:Bfamily}
 B=\diag(I_r,0),\qquad R=\begin{pmatrix}I_r\\0\end{pmatrix}.
\end{equation}
In particular, $\rank B=r$, and all the hypotheses of theorem \ref{thm:main} hold.

\begin{lemma}\label{lem:open}
For $G_r$, the matrices $X_s$ from \eqref{eq:PH}, with every coordinate
of $s$ strictly positive, have an image containing a nonempty open
subset of $\SPD_r$.
\end{lemma}

\begin{proof}
Write $s_i$ for the potential at $v_i$, and $s_{ij}$ for the potential
at $u_{ij}$. For $r=1$ there are no auxiliary vertices,
$H_s=(\kappa+s_1)^{-1}$, and varying $s_1>0$ gives an open
interval. Assume henceforth that $r\geq2$.

Keep the ordering with the interface vertices first.
Let $q=\binom r2$, and let $J\in\R^{r\times q}$ have column
\[
 J_{\cdot,(i,j)}=e_i+e_j
 \qquad(1\leq i<j\leq r).
\]
Thus $J$ records which interface vertices are joined to each
auxiliary vertex.

Each $v_i$ has $r-1$ auxiliary neighbors and one reflected
neighbor. Its diagonal entry in $A$ is therefore $r+\kappa$;
subtracting $B$ removes the contribution $1$ of the crossing
edge. Each $u_{ij}$ has degree two, and its row in $B$ is zero.
Consequently, if
\[
 D_s=\diag(r-1+\kappa+s_1,\ldots,r-1+\kappa+s_r),
 \qquad
 E_s=\diag(2+\kappa+s_{ij}:i<j),
\]
then
\[
 P_s=
 \begin{pmatrix}
 D_s&-J\\
 -J^{\T}&E_s
 \end{pmatrix},
 \qquad
 R=\begin{pmatrix}I_r\\0\end{pmatrix}.
\]
The Schur complement of the auxiliary block $E_s$ in $P_s$ is
\[
 M_s=D_s-JE_s^{-1}J^{\T}.
\]
It is the effective matrix on the interface obtained by eliminating
the auxiliary coordinates. We call $M_s$ the
\emph{interface Schur complement}.

The upper-left $r\times r$ block of $P_s^{-1}$ is
$H_s=R^{\T}P_s^{-1}R$.

For illustration, take $r=2$. The positive half has vertices
$v_1,v_2,u_{12}$, with $u_{12}$ joined to both $v_1$ and $v_2$.
Put
\[
 a_i=1+\kappa+s_i\quad(i=1,2),
 \qquad d=2+\kappa+s_{12}.
\]
In the order $(v_1,v_2,u_{12})$,
\[
 P_s=
 \begin{pmatrix}
 a_1&0&-1\\
 0&a_2&-1\\
 -1&-1&d
 \end{pmatrix},
 \qquad
 R=
 \begin{pmatrix}
 1&0\\
 0&1\\
 0&0
 \end{pmatrix}.
\]
To compute $H_sb=R^{\T}P_s^{-1}Rb$, solve
\[
 P_s
 \begin{pmatrix}x_1\\x_2\\y\end{pmatrix}
 =
 \begin{pmatrix}b_1\\b_2\\0\end{pmatrix}.
\]
The auxiliary equation is
\[
 -x_1-x_2+dy=0,
 \qquad\text{so}\qquad
 y=\frac{x_1+x_2}{d}.
\]
Substituting this into the first two equations gives
\[
 \begin{pmatrix}
 a_1-d^{-1}&-d^{-1}\\
 -d^{-1}&a_2-d^{-1}
 \end{pmatrix}
 \begin{pmatrix}x_1\\x_2\end{pmatrix}
 =
 \begin{pmatrix}b_1\\b_2\end{pmatrix}.
\]
Since $H_sb=(x_1,x_2)^{\T}$, this gives $H_s=M_s^{-1}$ with
\[
 M_s=\diag(a_1,a_2)
     -\frac{(e_1+e_2)(e_1+e_2)^{\T}}{2+\kappa+s_{12}}.
\]
This rank-two example already displays the full mechanism. Writing
$d=2+\kappa+s_{12}$, the three independent variations of $M_s$ are
\[
 \frac{\partial M_s}{\partial s_1}
   =\begin{pmatrix}1&0\\0&0\end{pmatrix},\qquad
 \frac{\partial M_s}{\partial s_2}
   =\begin{pmatrix}0&0\\0&1\end{pmatrix},\qquad
 \frac{\partial M_s}{\partial s_{12}}
   =\frac1{d^2}\begin{pmatrix}1&1\\1&1\end{pmatrix}.
\]
The first two control the diagonal entries, while the third supplies a
nonzero off-diagonal entry. They therefore span $\Sym_2$.

The same elimination works in general. To solve
$P_s\binom{x}{y}=\binom{b}{0}$, use the second block equation
to write $y=E_s^{-1}J^{\T}x$. The first then becomes
$(D_s-JE_s^{-1}J^{\T})x=b$. Thus $H_s=M_s^{-1}$, and the
columns of $J$ give
\begin{equation}\label{eq:Ms}
 M_s=\diag(r-1+\kappa+s_1,\ldots,r-1+\kappa+s_r)
      -\sum_{i<j}
       \frac{(e_i+e_j)(e_i+e_j)^{\T}}{2+\kappa+s_{ij}}.
\end{equation}
The Schur complement $M_s$ is positive definite since $P_s\succ0$.
Its derivatives are
\[
 \frac{\partial M_s}{\partial s_i}=e_ie_i^{\T},\qquad
 \frac{\partial M_s}{\partial s_{ij}}
 =\frac{(e_i+e_j)(e_i+e_j)^{\T}}{(2+\kappa+s_{ij})^2}.
\]

The derivatives with respect to $s_1,\ldots,s_r$ give the diagonal
basis matrices $e_i e_i^{\T}$. For each $i<j$, we also obtain the
corresponding off-diagonal basis matrix by subtracting the two
diagonal contributions:
\[
 (2+\kappa+s_{ij})^2\frac{\partial M_s}{\partial s_{ij}}
 -\frac{\partial M_s}{\partial s_i}
 -\frac{\partial M_s}{\partial s_j}
 =e_i e_j^{\T}+e_j e_i^{\T}.
\]
Thus the derivatives span $\Sym_r$. Since there are
$r+\binom{r}{2}=\dim\Sym_r$ of them, they form a basis.

The differential of matrix inversion is
\[
    \dd(M^{-1})[K]=-M^{-1}KM^{-1}.
\]
For $M\succ0$ this is an invertible linear map of $\Sym_r$ onto itself.
Thus the Jacobian of $s\mapsto X_s=I_r/2+M_s^{-1}$ also has full rank.
The inverse function theorem gives an open image near every strictly
positive $s$.
\end{proof}

\begin{proof}[Proof of Theorem~\ref{thm:main}\textup{(ii)}]
Sufficiency follows from part~(i). Conversely, suppose the
occupation field on $G_r$ is reflection positive.
Lemma~\ref{lem:laplace} makes the Laplace kernel positive definite.
In \eqref{eq:factor}, all factors other than
$\det(X_s+X_t)^{-\alpha}$ have the form $b(s)b(t)$ with $b>0$.
Conjugating each finite Gram matrix by the positive diagonal matrix
$\diag(b(s_i)^{-1})$ therefore shows that
$(s,t)\mapsto\det(X_s+X_t)^{-\alpha}$ is positive definite.
By Lemma~\ref{lem:open}, the image of $s\mapsto X_s$ contains
an open set $U\subset\SPD_r$. Choosing preimages for any finite
collection of points of $U$ shows that
$(X,Y)\mapsto\det(X+Y)^{-\alpha}$ is positive definite on $U$.
Proposition~\ref{prop:local} therefore gives $\alpha\in\W_r$.
\end{proof}

\begin{proof}[Proof of Theorem~\ref{thm:main}\textup{(iii)}]
Every nonnegative half-integer belongs to every $\W_r$, so
part~(i) gives sufficiency. If $\alpha>0$ is not a half-integer,
set $r=\lfloor2\alpha\rfloor+2$. Then $\alpha\notin\W_r$,
and part~(ii) provides a counterexample with unit conductances
and any prescribed constant positive killing.
For example, every $\alpha\in(1/2,1)$ is ruled out by the graph $G_3$,
while every $\alpha\in(1,3/2)$ is ruled out by $G_4$.
\end{proof}

\section{A band of finite width}\label{sec:bands}

Take the band $\mathbb R\times\{1,\ldots,r\}$ with reflection
$\rho(t,x)=(-t,x)$. The longitudinal motion has generator
$\partial_t^2$. Independently, the transverse coordinate jumps at
rate one to each neighbour and is killed on reaching $0$ or $r+1$.
There is no interior killing. The positive half is $t>0$.
Local times are occupation densities with respect to Lebesgue
measure in $t$ and counting measure in the transverse coordinate.

\begin{proposition}\label{prop:infinite-band}
On this band the loop-soup occupation field is reflection
positive if and only if $\alpha\in\W_r$.
\end{proposition}

The idea is to apply the determinant criterion of
Proposition~\ref{prop:local} to the joint Laplace transform of
finitely many occupations and their reflections. Each observation
point gives one rank-one variation of the matrix in that criterion.
We shall see that these variations span $\Sym_r(\R)$.

\begin{proof}
Put $J=r+1$. The transverse Dirichlet eigenfunctions $u_j$ and their
eigenvalues $\omega_j^2$ are given by
\begin{equation}\label{eq:band-dirichlet-spectrum}
 u_j(x)=\sqrt{2/J}\sin(\pi jx/J),\qquad
 \omega_j=2\sin(\pi j/(2J)),\qquad 1\le j\le r .
\end{equation}
With $a_j=(2\omega_j)^{-1/2}$, the Green function of
$-\partial_t^2+\omega_j^2$ gives
\begin{equation}\label{eq:infinite-band-green}
 G((t,x),(t',y))=\sum_{j=1}^r a_j^2 u_j(x)u_j(y)e^{-\omega_j|t-t'|}.
\end{equation}

For $z=(t,x)$ in the positive half put
$h(z)=(a_j u_j(x)e^{-\omega_jt})_{j=1}^r$.
Then $G(z,\rho w)=h(z)^{\T}h(w)$. Choose finitely many such
points $E$. Let $V$ have rows $h(z)^{\T}$, $z\in E$, and put
$D=(G(z,w)-G(z,\rho w))_{z,w\in E}$. The full Green matrix
on $E,\rho E$ is
\[
 \mathcal G=\diag(D,D)+
 \begin{pmatrix}V\\ V\end{pmatrix}
 \begin{pmatrix}V\\ V\end{pmatrix}^{\T}.
\]
Its quadratic form on vectors $(v,-v)$ gives $D\succeq0$.

Write $L_E=(L_\alpha(z))_{z\in E}$ and
$L_{\rho E}=(L_\alpha(\rho z))_{z\in E}$ for the point local times.
Their joint Laplace transform is $\det(I+\mathcal G\,\diag(S,T))^{-\alpha}$,
where $S=\diag s$, $T=\diag t$, and $s,t\geq0$.
On each layer the mean $\alpha G(z,z)$ and covariance
$\alpha G(z,w)^2$ are continuous. Thus local averages recover point
values in $L^2$, including the value at $t=0$ from the positive side.
Conversely, Riemann sums recover smooth compactly supported tests.
The matrix determinant lemma therefore gives
\begin{equation}\label{eq:band-laplace}
 \E_\alpha e^{-\langle s,L_E\rangle-\langle t,L_{\rho E}\rangle}
 =\frac{\det(X_s+X_t)^{-\alpha}}
        {\det(I+DS)^\alpha\det(I+DT)^\alpha},
 \qquad X_s=\tfrac12 I_r+V^{\T}S(I+DS)^{-1}V .
\end{equation}
Here $X_s\succ0$, since $S(I+DS)^{-1}=(S^{-1}+D)^{-1}$ for
$S\succ0$, with extension by continuity. For $\alpha\in\W_r$
the numerator is a positive definite kernel, and the denominator
has a separate positive factor in each argument. Sufficiency
follows by the density of the Laplace functions for each finite $E$, as in
Lemma~\ref{lem:laplace}, and then by the density of cylindrical
functions in $L^2$ of the positive-half field.

\smallskip\noindent\emph{Necessity.}
The response matrix $X_s$ is analytic near $s=0$. Since
\[
 S(I+DS)^{-1}=S-SDS+O(|s|^3),
\]
varying the potential at one observation point $z$ gives the first
variation
\begin{equation}\label{eq:band-first-variation}
 \left.\frac{\partial X_s}{\partial s_z}\right|_{s=0}
      =h(z)h(z)^{\T}.
\end{equation}
We need enough points for these rank-one matrices to span all of
$\Sym_r$. Equivalently, the scalar functions
\[
    z\longmapsto h_i(z)h_j(z),\qquad 1\leq i\leq j\leq r,
\]
must be linearly independent.

Suppose the functions $h_i h_j$ satisfy a linear relation. Absorbing
the positive constants $a_i a_j$ into the coefficients, we may write it as
\begin{equation}\label{eq:band-function-relation}
 \sum_{1\le i\le j\le r} c_{ij}u_i(x)u_j(x)e^{-(\omega_i+\omega_j)t}=0
 \qquad(t>0,\ 1\le x\le r).
\end{equation}
Summing over $x$ and using orthonormality gives
$\sum_i c_{ii}e^{-2\omega_i t}=0$. Since $\omega_1<\cdots<\omega_r$,
multiplying by $e^{2\omega_1t}$ and letting $t\to\infty$ gives $c_{11}=0$.
Repeating in increasing order gives $c_{ii}=0$ for every $i$.

We recover the successive diagonals in the same way. Put
$w(x)=2\cos(\pi x/J)$. The sine addition formula gives
\[
 w u_i=u_{i-1}+u_{i+1},\qquad u_0=u_{r+1}=0.
\]
Each multiplication by $w$ changes the index by one. Thus $w^d u_i$
only involves indices at most $i+d$. When $i+d\le r$, the coefficient
of $u_{i+d}$ is one, since reaching it requires $d$ successive increases.

Let $1\le d<r$ and suppose that $c_{ij}=0$ for $j-i<d$. Multiply
\eqref{eq:band-function-relation} by $w(x)^d$ and sum over $x$.
By orthonormality, $\sum_x w(x)^d u_i(x)u_j(x)$ is the coefficient
of $u_j$ in $w^d u_i$. Hence the relation reduces to
\[
 \sum_{i=1}^{r-d} c_{i,i+d}\,e^{-(\omega_i+\omega_{i+d})t}=0
 \qquad(t>0).
\]
For fixed $d$, the rates $\omega_i+\omega_{i+d}$ increase strictly
with $i$. Linear independence of these exponentials gives
$c_{i,i+d}=0$ for every $i$. Induction on $d$ proves the claim.

Let $q=r(r+1)/2$. Linear independence of the $q$ functions implies that
there are points $z_1,\ldots,z_q$ for which the matrices
$h(z_a)h(z_a)^{\T}$ are linearly independent. 
%Otherwise all evaluation
%vectors would lie in a proper subspace of $\R^q$, producing a nontrivial
%relation of the form \eqref{eq:band-function-relation}. 
By
\eqref{eq:band-first-variation}, the Jacobian of
\[
    (s_{z_1},\ldots,s_{z_q})\longmapsto X_s
\]
is invertible at $s=0$. The Jacobian determinant remains nonzero at a
nearby point with every coordinate strictly positive. The inverse
function theorem then shows that the image contains an open subset of
$\SPD_r$.

Finally, reflection positivity makes the left-hand side of
\eqref{eq:band-laplace} a positive definite kernel. The two denominator
factors depend separately and positively on $s$ and $t$, so they can be
removed by diagonal conjugation of each finite Gram matrix. Pulling the
remaining kernel through the local inverse just obtained shows that
$\det(X+Y)^{-\alpha}$ is positive definite on an open subset of
$\SPD_r$. Proposition~\ref{prop:local} gives $\alpha\in\W_r$.
\end{proof}

\begin{remark}
On the unit-conductance square-lattice band with half-integer columns,
the same proof applies. In \eqref{eq:infinite-band-green}, replace
the exponential rate $\omega_j$ and the prefactor $a_j^2$, respectively, by
\[
 \operatorname{arcosh}(1+\omega_j^2/2)
 =2\log\frac{\omega_j+\sqrt{4+\omega_j^2}}{2},
 \qquad \frac{1}{\omega_j\sqrt{4+\omega_j^2}} .
\]
The lattice rates still increase strictly with $j$, so the same
successive limits along the positive half-integer columns separate the
terms on each diagonal.
\end{remark}

In the next section we consider the continuous strip
$\mathbb R\times(0,\pi)$ and the renormalized occupation field of
Brownian motion killed at its boundary. The transverse coordinate
now ranges over the interval $(0,\pi)$, and every transverse sine
eigenfunction contributes to the reflected Green function. We show
how to recover the same quantity $D_XD_Y\det(X+Y)^{-\alpha}$ at $X=Y=I_m/2$
in each finite dimension $m$ by varying smooth potentials.

\section{Renormalized Brownian occupation fields}\label{sec:brownian}
%A scaling limit of occupation fields preserves reflection positivity
%for bounded continuous tests. To deduce failure of positivity by this
%route one would need convergence of specific negative reflection forms.
%We instead work directly with the renormalized Laplace transform.

\subsection{The centered field and its Laplace transform}

Let $D\subset\R^d$, $d=2$ or $3$, be a bounded domain, and let
\[
 H_D=-\tfrac12\Delta_D+m^2,\qquad G_D=H_D^{-1},\qquad m^2\geq0,
\]
with Dirichlet boundary conditions. We also use the massless planar
strip below. Write $\mu_D$ for the Brownian loop measure, including
the killing, and $\mathcal L_\alpha$ for its Poisson soup.
For a real smooth compactly supported test function $f$, the
centered occupation field is the $L^2$ limit
\begin{equation}\label{eq:brownian-centered}
 X_\alpha(f)=\lim_{\delta\downarrow0}
 \left[
 \sum_{\substack{\ell\in\mathcal L_\alpha\\T(\ell)>\delta}}
       \int_0^{T(\ell)}f(\ell_s)\,\dd s
 -\alpha\int_{T(\ell)>\delta}
       \int_0^{T(\ell)}f(\ell_s)\,\dd s\,\mu_D(\dd\ell)
 \right].
\end{equation}
The Poisson covariance formula and the second loop occupation moment
give the following finite integral, since $G_D(x,y)^2$ is locally
integrable in dimensions below four.
\begin{equation}\label{eq:brownian-covariance}
 \E[X_\alpha(f)X_\alpha(g)]
 =\alpha\int_{D\times D}f(x)G_D(x,y)^2g(y)\,\dd x\,\dd y.
\end{equation}
For $f\geq0$,
\begin{equation}\label{eq:brownian-laplace}
 \E e^{-X_\alpha(f)}
 =\det\nolimits_2(I+M_{\sqrt f}G_DM_{\sqrt f})^{-\alpha},
 \qquad
 \det\nolimits_2(I+K)=\prod_j(1+\lambda_j(K))e^{-\lambda_j(K)}.
\end{equation}
Here $M_h$ denotes multiplication by $h$. These statements and
the Gaussian identification
\begin{equation}\label{eq:brownian-isomorphism}
 X_{N/2}\ \overset{\mathrm{law}}=
       \frac12\sum_{a=1}^N:\!\phi_a^2\!:,
                 \qquad N\in\{1,2,\ldots\},
\end{equation}
are established in \cite[Section~10.2, Theorem~8 and
Corollary~7]{LeJan}; the $\phi_a$ are independent real Gaussian
fields with covariance $G_D$. For compactly supported tests on the strip,
one may exhaust $S$ by bounded Dirichlet rectangles. On each compact
observation set their Green kernels increase to $G_S$ and converge
in Hilbert--Schmidt norm by local square integrability of $G_S$.
Coupling the soups by restriction gives $L^2$ convergence of their
compensated occupations by the Poisson isometry. Applying
\eqref{eq:brownian-laplace} also to $2f$ gives uniform integrability
of the Laplace tests, so the identity passes to the limit.

Suppose $D$ and the killing are invariant under reflection $\rho$
in a hyperplane. Reflection positivity means
\begin{equation}\label{eq:brownian-rp}
 \E_\alpha[\overline{F(\rho X_\alpha)}F(X_\alpha)]\geq0
\end{equation}
for every bounded cylindrical functional with smooth test
functions supported in the positive half. Here
$(\rho X)(f)=X(f\circ\rho)$. These cylindrical functions generate the positive-half
$\sigma$-field. Density, reflection invariance and Cauchy--Schwarz
extend the inequality to its $L^2$ space.

On the symmetric domains used below, standard Gaussian reflection
positivity and \eqref{eq:brownian-isomorphism} give sufficiency for
$\alpha=N/2$ for $N\geq1$. At intensity zero the soup is empty, so $X_0=0$.

\subsection{The massless planar strip}\label{subsec:strip}

Let $S=\mathbb R\times(0,\pi)$, with generator $\tfrac12\Delta$,
Dirichlet boundary conditions and reflection $\rho(t,x)=(-t,x)$.

\begin{theorem}\label{thm:strip}
The centered Brownian occupation field on $S$ is reflection
positive if and only if $2\alpha$ is a nonnegative integer.
\end{theorem}

%As on the band of finite width, we vary a scalar potential in the
%positive half and apply the matrix gamma criterion to the reflected
%Laplace transform. On the strip the resulting operator has infinitely
%many transverse coordinates. What we need to recover, however, is
%only the derivative of order $m$ on each side of the determinant
%kernel on an $m\times m$ block. We show how to recover this derivative
%using smooth potentials supported in a single rectangle.
%
%The necessity argument has four steps:
%\begin{enumerate}[label=\textup{\arabic*.}]
%\item factor the reflected Laplace transform into positive one-sided
%factors and a determinant kernel depending on a response operator $R(f)$;
%\item show that the first variations of $R(f)$ generated by smooth scalar
%potentials are dense in the self-adjoint trace-class operators;
%\item correct the higher Taylor coefficients so that, to the required
%finite order, $R(f(Z))$ behaves like a prescribed $m\times m$ matrix $Z$;
%\item transfer differential positivity to
%$\det(I_m+Z+W)^{-\alpha}$ and use the negative Wallach polynomial.
%\end{enumerate}
%Only a finite Taylor jet is required at each fixed $m$.
%
\begin{proof}
The preceding subsection gives sufficiency for every nonnegative
half-integer intensity. For necessity assume reflection positivity and
fix an open rectangle $U$ with compact closure in the positive
half-strip. Put $e_j(x)=\sqrt{2/\pi}\sin(jx)$. Separation of
variables gives the off-diagonal Green kernel
\begin{equation}\label{eq:strip-Green}
 G_S((t,x),(s,y))=\sum_{j\ge1}\frac{e_j(x)e_j(y)}{j}e^{-j|t-s|}.
\end{equation}
For $z=(t,x)\in U$ set $h_j(z)=e_j(x)e^{-jt}/\sqrt j$ and
let $H:\ell^2(\mathbb N;\mathbb R)\to L^2(U;\mathbb R)$ have columns $h_j$.
%We use real
%Hilbert spaces for the response operators and trace duality. 
Then
$G_S(z,\rho w)=\langle h(z),h(w)\rangle_{\ell^2}$.
The distance of $U$ from the reflection line makes $H$
Hilbert--Schmidt. Write $A_U$ for the integral operator on
$L^2(U)$ with kernel $G_S(z,w)$ for $z,w\in U$. 
This operator is Hilbert--Schmidt by local square
integrability. Put $D=A_U-HH^*$. Identifying $\rho U$ with $U$, the Green
operator on $U\cup\rho U$ is
\[
 \mathcal G=\diag(D,D)+\mathcal H\mathcal H^*,
 \qquad \mathcal H=\binom H H .
\]
In particular $D\geq0$ by restriction to vectors $(v,-v)$.

%\smallskip\noindent\emph{The determinant in the reflection form.}
For a smooth potential $f$ supported in $U$, write $M_f$ for
multiplication by $f$. For $f\geq0$, or for real $f$ sufficiently
small in $L^\infty(U)$, set
\begin{equation}\label{eq:strip-response}
 R(f)=H^*M_f(I+DM_f)^{-1}H .
\end{equation}
For $f\geq0$ the middle operator equals
$M_{\sqrt f}(I+M_{\sqrt f}DM_{\sqrt f})^{-1}M_{\sqrt f}$,
so $R(f)$ is positive. For nonnegative $f,g$, formula
\eqref{eq:brownian-laplace} gives
\begin{equation}\label{eq:strip-laplace-factor}
 \E_\alpha e^{-X_\alpha(f)-X_\alpha(g\circ\rho)}
 =b(f)b(g)\det\bigl(I+R(f)+R(g)\bigr)^{-\alpha},
\end{equation}
where $b(f)>0$ depends only on $f$. Indeed, put $M=\diag(M_f,M_g)$.
Cycling the factors in $\det{}_2(I+M^{1/2}\mathcal G M^{1/2})$
gives $\det{}_2(I+\mathcal G M)$. Write $\mathcal G M=A+B$ with
$A=\diag(D,D)M$ and $B=\mathcal H\mathcal H^*M$. Here $A$ is Hilbert--Schmidt
and $B$ is trace class. The identity
\[
 \det{}_2(I+A+B)=
 \det{}_2(I+A)\,e^{-\tr B}\,
 \det\bigl(I+(I+A)^{-1}B\bigr)
\]
follows from the finite matrix identity by approximation in the
respective norms. The first two factors split over the two halves.
In particular
$b(f)=\det{}_2(I+DM_f)^{-\alpha}
       \exp\{\alpha\tr(H^*M_fH)\}$.
In the last determinant, interchanging the factors moves it to
$\ell^2$ and gives $I+R(f)+R(g)$. This proves
\eqref{eq:strip-laplace-factor}. The Laplace tests are in $L^2$
by \eqref{eq:brownian-laplace} applied to $2f$. Reflection
positivity and an easy truncation argument therefore make the determinant kernel
in \eqref{eq:strip-laplace-factor} positive definite on
nonnegative smooth potentials.

\smallskip\noindent\emph{The first variations.}
Near zero, $R$ is analytic from real bounded potentials to real
self-adjoint trace-class operators, by the Neumann series in
\eqref{eq:strip-response}. Its derivative is
\[
 L(f)=H^*M_fH=\int_U f(z)\,h(z)h(z)^*\,\dd z .
\]
The range of $L$ on real $C_c^\infty(U)$ is dense in this
trace-class space for the trace norm. To see this, let a bounded
self-adjoint operator $T$ annihilate the range under the trace
pairing. Testing against every $f$ and using continuity of
$h:U\to\ell^2$ gives $\langle Th(z),h(z)\rangle=0$ throughout $U$.
For fixed $x$ in the transverse interval the series for this scalar
is a power series in $u=e^{-t}$, absolutely convergent for $|u|<1$
since $|T_{ij}|\leq\|T\|$. It vanishes on an interval of $u$,
so its coefficient of each degree $k$ vanishes.
\[
 \sum_{\substack{i+j=k\\i\le j}}(2-\delta_{ij})
       \frac{T_{ij}}{\sqrt{ij}}\sin(ix)\sin(jx)=0 .
\]
This finite trigonometric polynomial vanishes on an interval.
For the unordered pairs $i\le j$ with $i+j=k$, the frequencies
$j-i$ are distinct and strictly less than $k$. The identity
$2\sin(ix)\sin(jx)=\cos((j-i)x)-\cos(kx)$ therefore gives
$T_{ij}=0$ for every pair: each frequency $j-i<k$ occurs in only
one summand, while the frequency $k$ is the common second term.
Trace-class duality proves density.

For a concrete finite-dimensional picture, take $m=2$ and write
\[
 Z=z_{11}E_{11}+z_{12}(E_{12}+E_{21})+z_{22}E_{22}.
\]
Density allows us to choose three smooth potentials whose first
variations under $L$ approximate the three displayed matrix directions
in trace norm. Thus scalar potentials already generate, to first order,
all three independent entries of a symmetric $2\times2$ matrix. The
next step corrects the nonlinear terms of $R$ without changing the
coefficients already fixed at lower degree.

\smallskip\noindent\emph{Recovering the finite matrix derivative.}
Fix $m$ and embed a real symmetric $m\times m$ matrix $Z$ as
the leading block $\iota Z$ on $\ell^2$, zero elsewhere.
For every $\varepsilon>0$ we can choose a polynomial $f(Z)$
in the independent entries of $Z$, with coefficients in
$C_c^\infty(U)$ and $f(0)=0$, such that every Taylor coefficient
of $R(f(Z))-\iota Z$ through degree $m$ has trace norm less than
$\varepsilon$.

Here is the recursive construction in more detail. Write
\[
    f(Z)=\sum_{1\leq|\nu|\leq m}Z^\nu f_\nu,
    \qquad f_\nu\in C_c^\infty(U),
\]
where $\nu$ ranges over monomials in the $m(m+1)/2$ independent
entries of $Z$. The Neumann series gives
\begin{equation}\label{eq:strip-response-series}
 R(f)=L(f)-H^*M_fDM_fH
      +H^*M_fDM_fDM_fH-\cdots .
\end{equation}
At degree one, density of $L$ lets us choose the $f_\nu$ so that the
linear coefficient approximates the corresponding coefficient of
$\iota Z$. Suppose all coefficients of degree less than $q$ have been
chosen. For a monomial $Z^\nu$ of degree $q$, its coefficient in
$R(f(Z))$ has the form
\[
    L(f_\nu)+A_\nu,
\]
where $A_\nu$ is completely determined by the already chosen
lower-degree coefficients through the nonlinear terms in
\eqref{eq:strip-response-series}. Choose $f_\nu$ so that
$L(f_\nu)$ approximates the desired coefficient minus $A_\nu$.
Proceeding degree by degree through $q=m$ produces the required jet.
For example, along $f(s)=sf_1+s^2f_2$, the quadratic coefficient is
\[
    L(f_2)-H^*M_{f_1}DM_{f_1}H,
\]
so we choose $f_2$ to make this difference smaller than $\varepsilon$
in trace norm. Only finitely many
monomials occur, and the approximation tolerances may be divided among
them so that every coefficient error is below $\varepsilon$.
Let $D_Z$ be the determinant differential operator of
Lemma~\ref{lem:gamma-differential}. 
The reason that the construction stops at degree $m$ is that $D_Z$ has
total differential order $m$. Consequently
$D_ZD_W$ at $Z=W=0$ depends only on the $m$-jet in each variable.

For each such polynomial,
\begin{equation}\label{eq:strip-differential}
 \left.D_ZD_W\det\bigl(I+R(f(Z))+R(f(W))\bigr)^{-\alpha}
       \right|_{Z=W=0}\geq0 .
\end{equation}
Let $K$ be the union of the supports of the coefficients of $f$, and
choose a nonnegative $\chi\in C_c^\infty(U)$ equal to one on $K$.
Fix $\delta>0$. Since $f$ is a polynomial with $f(0)=0$, we have
$\|f(Z)\|_\infty<\delta/2$ on some neighbourhood of $Z=0$.
There the potential $g_\delta(Z)=\delta\chi+f(Z)$ is at least
$\delta/2$ on $K$ and equals $\delta\chi\geq0$ outside $K$.
Hence \eqref{eq:strip-laplace-factor} gives a positive definite
kernel on this neighbourhood,
\[
 k_\delta(Z,W)
 =\det\bigl(I+R(g_\delta(Z))+R(g_\delta(W))\bigr)^{-\alpha}.
\]

By analyticity of $R$ and of the Fredholm determinant, this expression
is jointly analytic in $(\delta,Z,W)$ near $(0,0,0)$, where the
determinant equals one. For each sufficiently small fixed $\delta>0$,
\eqref{eq:differential-positivity} therefore gives
$(D_ZD_W k_\delta)(0,0)\geq0$. This derivative is continuous in
$\delta$ at zero, so letting $\delta\downarrow0$ gives
\eqref{eq:strip-differential}.

%To justify the sign, choose a nonnegative $\chi\in C_c^\infty(U)$
%equal to one on the supports of all coefficients of $f$.
%For each $\delta>0$ small enough, $\delta\chi+f(Z)$ is
%nonnegative for $Z$ in a neighbourhood of zero.
%Apply \eqref{eq:differential-positivity} there to the kernel
%in \eqref{eq:strip-laplace-factor}, with this potential on
%each side. For the fixed polynomial $f$, analyticity near zero lets
%$\delta\downarrow0$ after differentiation and gives \eqref{eq:strip-differential}.

The derivative in \eqref{eq:strip-differential} depends
continuously on the Taylor coefficients through degree $m$ on
each side, since the Fredholm determinant is analytic in the
trace norm. Letting $\varepsilon\downarrow0$ and applying the
matrix gamma calculation now gives
\[
 0\leq
 \left.D_ZD_W\det(I_m+Z+W)^{-\alpha}\right|_{Z=W=0}
 =p_m(\alpha)p_m(\alpha+1).
\]
If $2\alpha$ is not an integer, take $m=\lfloor2\alpha\rfloor+2$.
Then $p_m(\alpha)$ has exactly one negative factor and
$p_m(\alpha+1)>0$, a contradiction. For example, the $m=2$
calculation excludes $0<\alpha<1/2$, and the $m=3$ calculation excludes
$1/2<\alpha<1$. Higher values of $m$ exclude the successive open
intervals between nonnegative half-integers.
\end{proof}

\section{Extensions to bounded reflection-invariant domains}
\label{sec:bounded-domains}

We outline how the necessity of the half-integer condition extends
from the strip to bounded connected reflection-invariant domains
with smooth boundary. In two dimensions we localize to a disk and
apply conformal covariance to transfer Theorem~\ref{thm:strip} to
the disk. In three dimensions we localize to a thin cylinder and
recover the occupation field on its planar base by dimension reduction.

\subsection{Symmetric changes of killing and localization}
\label{subsec:brownian-tilts}

On a bounded domain let $H=-\Delta_D/2+m^2$ and $G=H^{-1}$.
For bounded real $V$, $X_\alpha^H(V)$ is defined by approximation
in the covariance norm \eqref{eq:brownian-covariance}.
If $V$ is reflection invariant and
$H+M_V>0$, define the tilt
\begin{equation}\label{eq:brownian-tilt}
 \frac{\dd\mathbb P_{H,V}}{\dd\mathbb P_H}
 =\frac{e^{-X_\alpha^H(V)}}{\E_H e^{-X_\alpha^H(V)}} .
\end{equation}
Its normalizer is
$\det_2(I+G^{1/2}M_VG^{1/2})^{-\alpha}$. Under this tilt, the field is,
up to a deterministic change of centering, the occupation field for
$H+M_V$:
\begin{equation}\label{eq:brownian-recenter}
 X_\alpha^H(f)-c_V(f)\ \overset{\rm law}{=} X_\alpha^{H+M_V}(f),\qquad
 c_V(f)=\alpha\tr\big(((H+M_V)^{-1}-H^{-1})M_f\big),
\end{equation}
with the tilted law on the left.

We shall use without further proof the standard stability facts that a
reflection-invariant tilt of this form preserves reflection positivity,
and that the killing potentials
\[
    V_M=M1_{D\setminus U},\qquad M\longrightarrow\infty,
\]
give the Dirichlet field on a reflection-invariant Lipschitz
subdomain $U\Subset D$ in the limit, for tests supported in $U$.
Reflection positivity passes to this limit. On $U$, the original
constant mass may then be removed by the stable tilt $V=-m^2$. These are
the only consequences of tilting and additional killing needed below.

\subsection{Other planar domains}

\begin{proposition}\label{prop:planar-domains}
Let $D\subset\R^2$ be a bounded connected reflection-invariant
domain with smooth boundary. For Dirichlet Brownian motion with
constant killing $m^2\geq0$, the centered occupation field is
reflection positive if and only if $2\alpha\in\{0,1,2,\ldots\}$.
\end{proposition}

\begin{proof}[Proof outline]
The domain contains a small disk centered on the reflecting line,
by connectedness and reflection invariance. Section~\ref{subsec:brownian-tilts}
transfers reflection positivity to its massless Dirichlet field.
A similarity reduces to the unit disk $\mathbb D$. The conformal map
\[
 \Phi(z)=\frac{e^z-i}{e^z+i},\qquad
 0<\operatorname{Im}z<\pi,\qquad
 \Phi(-\overline z)=-\overline{\Phi(z)}
\]
takes the strip onto $\mathbb D$ and intertwines their reflections.
Conformal covariance of the compensated occupation field
\cite[Section~10.2, Proposition~31]{LeJan} gives, jointly in $f$,
\begin{equation}\label{eq:brownian-conformal}
 X_\alpha^{\mathbb D}(f)\ \overset{\rm law}{=} X_\alpha^S\big(|\Phi'|^2(f\circ\Phi)\big).
\end{equation}
This is a bijective correspondence of compactly supported smooth
tests on the respective positive halves and it intertwines reflection.
To make the necessity implication explicit, suppose that $2\alpha$ is
not an integer. The proof of Theorem~\ref{thm:strip} produces finitely
many positive-half strip tests and coefficients for which the reflection
form is negative. Pulling those tests through $\Phi$ gives the same
negative reflection form on the massless disk. This contradicts the
reflection positivity inherited from $D$ by localization and removal of
the mass. Hence $2\alpha$ must be an integer. Sufficiency is the standard
Gaussian consequence of \eqref{eq:brownian-isomorphism}.
\end{proof}

\subsection{Three dimensions}

\begin{proposition}\label{prop:three-dimensional}
Let $D\subset\R^3$ be a bounded connected reflection-invariant domain
with smooth boundary. For Dirichlet Brownian motion with constant
killing $m^2\geq0$, its centered occupation field is reflection positive if and
only if $2\alpha\in\{0,1,2,\ldots\}$.
\end{proposition}

\begin{proof}[proof outline by dimension reduction]
Localize to a thin product domain $U\times I_\varepsilon\Subset D$.
Here $U$ is a fixed planar disk centered on the reflecting plane,
$I_\varepsilon=(-\varepsilon/2,\varepsilon/2)$ is in a direction
tangent to that plane, and reflection acts in $U$. The covariance
singularity $|x-y|^{-2}$ on a planar slice is not locally
integrable there. We use volume tests constant in the thin direction.
The transverse
Dirichlet eigenvalues are $\lambda_n(\varepsilon)=\pi^2n^2/(2\varepsilon^2)$.
After removing the original mass, tilt by $V=-\lambda_1(\varepsilon)$.
This is stable because
\begin{equation}\label{eq:thin-cylinder-operator}
 H_\varepsilon
 =H_U\otimes I+I\otimes
   \big(-\tfrac12\partial_z^2-\lambda_1(\varepsilon)\big)>0,
 \qquad H_U=-\tfrac12\Delta_U .
\end{equation}
Multiplication by $f\otimes1$ preserves the transverse
eigenspaces. The regularized determinant therefore factors over
them, giving jointly for every finite family of planar tests
\begin{equation}\label{eq:thin-cylinder-decomposition}
 X_\alpha^{H_\varepsilon}(f\otimes1)
 \ \overset{\rm law}{=} \sum_{n\geq1}X_{\alpha,n}^{H_U+\mu_n}(f),\qquad
 \mu_n=\frac{\pi^2(n^2-1)}{2\varepsilon^2}.
\end{equation}
The fields on the right are independent, and the $n=1$ summand
is the massless planar field. Indeed, factorization holds for the
characteristic function of every real linear combination of the
tests. The variance of the sum over $n\ge2$ is bounded by
$\alpha\|f\|_\infty^2\sum_{n\ge2}\tr(H_U+\mu_n)^{-2}$.
The Dirichlet heat kernel is bounded by the free planar heat kernel, so
\[
 \tr(e^{-tH_U})=\int_U p_U(t,x,x)\,\dd x
 \leq\frac{|U|}{2\pi t}.
\]
Thus for $\mu>0$,
\[
 \tr\big((H_U+\mu)^{-2}\big)
 =\int_0^\infty t e^{-\mu t}\tr(e^{-tH_U})\,\dd t
 \leq\frac{|U|}{2\pi\mu}.
\]
Since $\sum_{n\ge2}(n^2-1)^{-1}=3/4$, the series converges in
$L^2$ and satisfies
\begin{equation}\label{eq:thin-cylinder-error}
 \E\left|\sum_{n\geq2}X_{\alpha,n}^{H_U+\mu_n}(f)\right|^2
 \leq\frac{3\alpha|U|\|f\|_\infty^2}{4\pi^3}\,\varepsilon^2 .
\end{equation}
Cutoffs at the ends of $I_\varepsilon$ give $f\otimes1$ by
approximation in the covariance norm. The bound applies jointly
to any fixed finite family of tests and their reflections, by
applying it to each component. Reflection positivity thus passes
through the tilt and the thin-cylinder limit to the massless
planar disk. Proposition~\ref{prop:planar-domains} gives necessity,
and the Gaussian argument following \eqref{eq:brownian-rp}
gives sufficiency.
\end{proof}

The author acknowledges the help of artificial intelligence.

%\begin{remark}
%Components exchanged in pairs by the reflection carry independent
%soups. Necessity applies whenever a reflection-invariant component
%meets the reflecting hyperplane. In whole space take concentric
%balls centered on that hyperplane and impose killing at rate $M$
%in the annulus between them. As $M\to\infty$ the inner ball is
%separated from the exterior and its Green kernel acquires the
%Dirichlet boundary condition. The tilt and the convergence argument
%use only the Green operator restricted to a compact set containing
%the annulus and the observation supports. Its kernel is square
%integrable there for positive mass in dimensions two and three and
%for zero mass in dimension three. In particular the loop integral
%of $Y^2$ is finite for each compactly supported annular potential,
%and $0\le e^{-Y}-1+Y\le Y^2/2$ justifies its nonnegative tilt.
%Necessity in these cases follows from the ball result. For
%sufficiency one may exhaust whole space by concentric symmetric
%Dirichlet balls and use local Green convergence on compact tests.
%Massless Brownian motion in the whole plane requires a
%separate prescription for large loops.
%\end{remark}
%

\end{document}